\documentclass[11pt,A4]{article}
\usepackage{amsmath,amsthm,amsfonts,amssymb,amscd, amsxtra,mathrsfs, color}
\usepackage[english]{babel}

\usepackage[margin=1.0in]{geometry}
\newtheorem{theorem}{Theorem}
\newtheorem{lemma}{Lemma}
\newtheorem{corollary}{Corollary}

\newtheorem{remark}{Remark}

\newcommand{\m}{\mathbb{M}}
\newcommand{\tpm}{T_{p}\mathbb{M}}

\begin{document}

\title{On the superlinear convergence of  Newton's method \\  on Riemannian manifolds}

\author{
Fernandes, T. A.  \thanks{Universidade Estadual do Sudoeste da Bahia, BA 45083-900, BR (e-mail:{\tt telesfernandes@uesb.edu.br}). The author was supported in part by UESB.} 
\and
Ferreira, O. P.  \thanks{ Universidade Federal de Goias, Goiania, GO 74001-970, BR (e-mail:{\tt orizon@ufg.br}). The author was supported in part by CNPq Grants 305158/2014-7 and FAPEG/GO.} 
\and
Yun, Y. J.  \thanks{Universidade Federal do Paran\'a, Curitiba, PR 81531-980,  BR (e-mail:{\tt jin@ufpr.br}). The author was supported in part by CNPq and CAPES.}
}
\maketitle
\begin{abstract}
In  this paper  we study the  Newton's method for  finding a singularity  of a differentiable vector field defined on a Riemannian manifold.  Under the assumption of invertibility of covariant derivative of the vector field at its singularity, we establish  the  well definition of the method  in a suitable neighborhood  of this  singularity. Moreover, we also show that  the generated sequence  by Newton method    converges for the solution with superlinear rate.  



\end{abstract}

\section{Introduction}
Applications of the concepts of Riemannian geometry in optimization  arise  when  optimization problems  are formulated as  problems of  finding a minimizers  of a real-valued  functions  defined on  smooth nonlinear manifolds.  Indeed, many optimization problems are naturally posed on Riemannian manifolds, which have a specific underlying geometric and algebraic structure that can be exploited to greatly reduce the computational cost  of obtaining the  minimizer.  For instance, to take advantage of the Riemannian geometric structure, it is preferable to treat certain constrained optimization problems as problems for finding singularities of gradient vector fields on Riemannian manifolds rather than using Lagrange multipliers or projection methods; see \cite{Smith1994,Absil2009}.  Accordingly, constrained optimization problems are regarded as unconstrained ones from the viewpoint of Riemannian geometry.  The early works dealing with this issue include  \cite{Luenberger1972,Gabay1982,Rapcsak1989,Udriste1994}. In the recent years  there has been increasing  interest in the development of geometric optimization algorithms which  exploit the differential structure of the nonlinear manifold, papers published on this topic include, but are not limited to \cite{Absil2014,Huang2015,LiLopezMartin-Marquez2009,LiWang2008,WangLiLopezYao2015,WangLiWangYao2015,HuLiYang2016,Ring2012,Argyros2014,Huang2015,Manton2015,Bittencourt2016,De2016}.  In this paper, instead of focusing on  problems of finding singularities of gradient vector fields on Riemannian manifolds, which includes finding local minimizers, we consider the more general problem of finding singularities of vector fields.  

It is well known that Newton's method is powerful tools for finding zeros of nonlinear functions in Banach spaces. Moreover, Newton's method also serves as a powerful theoretical tool with a wide range of applications in pure and applied mathematics; see, for example,  \cite{Adler2002,Nash1956,Moser1961}.  For this reasons, it has inspired several studies that deal with the issue of generalizing it  from the linear  context to the Riemannian setting, one could refer to related works, including \cite{Owren2000,Li2009,Absil2014,Argyros2009,Dedieu2000,LiWang2008,Manton2015,Miller2005,Schulz2014,Wang2011,Wang2009,Wang2012,Zhang2010,Smith1994,Udriste1994,Absil2009,Ferreira2012,Adler2002,Ferreira2002,Dedieu2003,Jiang2007,Li2006}. It is worth mentioning that, in all these  previous papers,   the analysis presented  on the Riemannian version of Newton's  method   has been used Lipschitz condition or Lipschitz-like conditions on the covariant derivative of the vector field.  In fact, all these papers are concerned to establish the quadratic rate of convergence of the method. It seems that  some kind of control on the covariant  derivative of the vector  field,  in a suitable neighborhood of its singularity,  is needed   for obtaining the quadratic convergence rate of the sequence generated by Newton's Method.   Hence the convergence analysis of the Newton method for finding singularity of vector field in Riemannian manifolds under  Lipschitz condition or Lipschitz-like conditions  are  well known. However,  we also know that in the linear context, whenever the derivative of the function that define the equation is nonsingular at the solution, the Newton's method has local convergence with superlinear rate, \cite[chapter 8, Theorem 8.1.10, p. 148]{Ortega1990}.  As far as we know, the local superlinear convergence analysis of Newton's method under mild assumption, namely, only invertibility  of the covariant derivative of the vector field at its singularity, is a new contribution of this paper. 

The organization of the paper is as follows. In Section 2, some notation and basic result used in the paper are presented. In Section 3 we provide  the local superlinear convergence analysis of Newton's method. Finally, Section 4  concludes the paper.
\section{Basic definition and auxiliary results}
In this section we recall some notations, definitions and basic properties of Riemannian manifolds used throughout the paper. They can be found in many introductory books on Riemannian Geometry, for example, in \cite{doCarmo1992} and \cite{Sakai1996}. 

Let $\mathbb{M}$ be a smooth manifold, denote the {\it tangent space} of $\m$ at $p$ by $\tpm$ and the {\it tangent bundle} of $\m$ by $T\m=\bigcup_{p\in\m}\tpm$.  The corresponding norm associated to the Riemannian metric $\langle \cdot ~, ~ \cdot \rangle$ is denoted by $\|  ~\cdot~ \|$. The Riemannian  distance  between $p$ and $q$   in a finite dimensional Riemannian manifold $M$ is denoted  by $d(p,q)$,  which induces the original topology on $\mathbb{M}$, namely,  $(\mathbb{M}, d)$ is a complete metric space.   The open ball of radius $r>0$ centred at $p$ is defined as  $B_{r}(p):=\left\lbrace q\in\m:d(p,q)<r\right\rbrace$.  Let  $\Omega \subset  \m$ be an open set and denote by ${\cal X}(\Omega)$ the space of  $C^1$ vector fields on $\Omega$.  Let $\nabla$ be the Levi-Civita connection associated to $(\mathbb{M}, \langle \cdot ~, ~ \cdot \rangle)$.  The covariant derivative of $X \in {\cal X}(\Omega)$ determined by $\nabla$ defines at each $p\in \Omega$ a linear map $\nabla X(p):\tpm\to\tpm$ given by $\nabla X(p)v:=\nabla_{Y}X(p),$
where $Y$ is a vector field such that $Y(p)=v$.    The {\it norm of a  linear map} $A:T_p \mathbb{M} \to T_p \mathbb{M}$  is defined by $\|A\|:=\sup \left\{ \|Av \|~:~ v\in T_p \mathbb{M}, \,\|v\|=1 \right\}$.  A vector field $V$ along a differentiable curve $\gamma$ in $M$ is said to be {\it parallel} if and only if $\nabla_{\gamma^{\prime}} V=0$. If $\gamma^{\prime}$ itself is parallel we say that $\gamma$ is a {\it geodesic}.  The restriction of a geodesic to a  closed bounded interval is called a {\it geodesic segment}. A geodesic segment joining $p$ to $q$ in $\mathbb{M}$ is said to be {\it minimal} if its length is equal to $d(p,q)$. If there exists a  unique geodesic segment  joining $p$ to $q$, then we denote it  by $\gamma_{pq}$.  For each $t \in [a,b]$, $\nabla$ induces an isometry, relative to $ \langle \cdot , \cdot \rangle  $, $P_{\gamma,a,t} \colon T _{\gamma(a)} {\mathbb{M}} \to T _{\gamma(t)} {\mathbb{M}}$ defined by $ P_{\gamma,a,t}\, v = V(t)$, where $V$ is the unique vector field on $\gamma$ such that
$$
 \nabla_{\gamma'(t)}V(t) = 0, \qquad V(a)=v, 
$$
the so-called {\it parallel transport} along  the geodesic segment   $\gamma$ joining  $\gamma(a)$ to $\gamma(t)$. Note also that $ P_{\gamma,\,b_{1},\,b_{2}}\circ P_{\gamma,\,a,\,b_{1}}=P_{\gamma,\,a,\,b_{2}}$ and  $P_{\gamma,\,b,\,a}=P^{-1}_{\gamma,\,a,\,b}$.  When there is no confusion we will  consider the notation $P_{pq}$  instead of $P_{\gamma,\,a,\,b}$ in the case when $\gamma$ is the unique geodesic segment joining $p$. A Riemannian manifold is {\it complete} if the geodesics are defined for any values of $t\in \mathbb{R}$. Hopf-Rinow's theorem asserts that any pair of points in a  complete Riemannian  manifold $\mathbb{M}$ can be joined by a (not necessarily unique) minimal geodesic segment.  Due to  the completeness of the Riemannian manifold $\mathbb{M}$, the {\it exponential map} $\exp_{p}:T_{p}\mathbb{M} \to \mathbb{M} $ can be  given by $\exp_{p}v\,=\, \gamma(1)$, for each $p\in \mathbb{M}$, where $\gamma$ is the geodesic defined by its position $p$ and velocity $v$ at $p$. Let $p\in \mathbb{M}$, the {\it injectivity radius} of $\mathbb{M}$ at $p$ is defined by
$$
i_{p}:=\sup\left\{ r>0~:~{\exp_{p}}_{\lvert_{B_{r}(o_{p})}} \mbox{ is\, a\, diffeomorphism} \right\},
$$
where $o_{p}$ denotes the origin of $T_{p}\mathbb{M}$ and $B_r(0_{p}):=\left\lbrace  v\in T_{p}\mathbb{M}:\parallel v-0_{p}\parallel <r\right\rbrace$, is called neighbourhood of injectivity of $p$. 
\begin{remark}\label{unicidadedageodesica}
Let $\bar{p}\in \mathbb{M}$. The above definition implies  that if $0<\delta<r_{\bar{p}}$, then $\exp_{\bar p}B_{\delta}(0_{ \bar p})=B_{\delta}( \bar p)$. Moreover, for all $p\in B_{\delta}(\bar p)$,  there exists a unique geodesic segment $\gamma$ joining  $p$ to $\bar p$, which is given by $\gamma_{p\bar p}(t)=\exp_{p}(t \exp^{-1}_{p} {\bar p})$, for all $t\in [0, 1]$.
\end{remark}

Next, we present the number $K_{p}$,  introduced in    \cite{Dedieu2003},  which measures  how fast the geodesics spread apart in $\mathbb{M}$. Let  $i_{p}$ for $p\in\mathbb{M}$ be the radius of injectivity of $\mathbb{M}$ at $p$. Consider the amount given by 
$$
K_{p}:=\sup\left \{  \dfrac{d(\exp_{q}u, \exp_{q}v)}{\parallel u-v\parallel} ~:~ q\in B_{i_{p}}(p), ~ u,\,v\in T_{q}\mathbb{M}, ~u\neq v, ~\parallel v\parallel\leq i_{p},
\parallel u-v\parallel\leq i_{p}\right\}.
$$
\begin{remark} In particular, when $u=0$ or more generally when $u$ and $v$ are on the same line through $0$, $d(\exp_{q}u, \exp_{q}v)=\parallel u-v\parallel$. Hence,  $K_{p}\geq1$,  for all $p\in\m$. Moreover, when $\m$ has non-negative sectional curvature, the geodesics spread apart less than the rays \cite[chapter 5]{doCarmo1992}, i.e., $d(\exp_{p}u, \exp_{p}v)\leq\| u-v\|$ and, in this case, $K_{p}=1$  for all $p\in\m$.
\end{remark}
Let  $X\in {\cal X}(\Omega)$   and $  \bar{p}\in \Omega$. Assume that  $0<\delta<r_{\bar{p}}$. Since  $exp_{\bar{p}}B_{\delta}(0_{\bar{p}})=B_{\delta}(\bar{p})$ and   there exists a unique geodesic joining  each $p\in B_{\delta}(\bar{p})$ to $\bar{p}$. Then  using \cite [ equality 2.3]{Ferreira2002}, for each $p\in   B_{\delta}(\bar{p})$ we obtain
\begin{equation}\label{diferenciabilidade}
X(p)= P_{\bar{p},p}X(\bar{p})+P_{\bar{p},p}\nabla X(\bar{p})\exp^{-1}_{\bar{p}}p+ d(p, \bar{p})r(p),  \quad \underset{p\to \bar{p}}{\lim}r(p)=0.
\end{equation}
We end this section with a well-known Banach's Lemma.
\begin{lemma} \label{le:Banach}
Let $B$ be a linear operator and $I_{p}$ the identity operator in $\tpm$. If $\parallel B-I_{p}\parallel<1$, then $B$ is invertible and $\parallel B^{-1}\parallel\leq\dfrac{1}{\parallel B-I_{p}\parallel}.$
\end{lemma}
{\it Throughout the paper $\mathbb{M}$ is a complete  Riemannian manifold of finite dimension.}
\section{Superlinear convergence of   Newton's method}
In this section  we  study  the  Newton's method to find a point $p\in \Omega$ satisfying the equation
\begin{equation} \label{eq:TheProblem}
 X(p)=0, 
\end{equation}
where  $X:\Omega\to T\mathbb{M}$ a differentiable vector field, and $\Omega\subset\m$ is an open set. The Newton's method to solve \eqref{eq:TheProblem} formally generates a sequence, with an initial point $p_0\in \Omega$, as follows.
\begin{equation} \label{eq:NM}
p_{k+1}=\exp_{p_{k}}(-\nabla X(p_{k})^{-1}X(p_{k})),  \qquad k=0, 1, \ldots. 
\end{equation}
{\it From now on, we assume that $p_*\in \Omega$ is a solution of \eqref{eq:TheProblem}}.  Our aim is to prove that,  under  the assumption  of nonsingularity of covariant derivative at the solution $p_*$, the  iterates  \eqref{eq:NM} starting in a suitable neighbourhoods of  $p_*$  are well defined and converges superlinearly to $p_*$.  To obtain  this result,  we begin by stating  an  important property of the parallel transport to our context.  It is worth to mention that,  to ensure  this property  we  use the  same ideas given in   the proof of \cite [Lemma 2.4, item (iv)]{LiLopezMartin-Marquez2009} for Hadamard manifolds with some minor  necessary technical adjustments to fit  on  any Riemannian manifold.  
\begin{lemma}\label{le:ContTranpPar} 
Let  $\bar{p} \in  \mathbb{M}$,  $0<\delta<i_{\bar{p}}$ and  $u\in T_{\bar{p}}\mathbb{M}$.  Then,  the vector field $F:B_{\delta}(\bar{p}) \to T\mathbb{M}$ defined by $F(p):=P_{\bar{p}p}u$ is continuous.
\end{lemma}
{\it Proof} Assume that $\mathbb{M}$  is $n$-dimensional.  Let $p\in B_{\delta}(\bar{p})$ and the unique geodesic segment $\gamma_{p}$, joining $\bar p$ to $p$ , given by remark \ref{unicidadedageodesica}. Let $u\in T_{\bar p}\mathbb{M}$. From  definition of the  parallel transport, there is unique continuously differentiable  vector filed $Y_{p}$ along $\gamma_{p}$ such that $Y_{p}\left(\gamma_{p}(0)\right)=u$, $Y_{p}\left(\gamma_{p}(1)\right)=P_{\bar p p}u$ and
\begin{equation}\label{derivadanula}
\nabla_{\gamma_{p}^{'}(t)}Y_{p}\left(\gamma_{p}(t)\right)=0,\qquad \forall ~ t\in[0,1], 
\end{equation}
see \cite[pag. 29]{Sakai1996}. The definition $i_{\bar p}$ implies that   $\varphi:=exp^{-1}_{\bar{p}}:B_{\delta}(\bar{p})\to B_{r_{\bar{p}}}(0_{\bar{p}})$ is a diffeomorphism and then  $\left(B_{\delta}(\bar{p}),\varphi\right)$ is a local chart at $\bar{p}$. For each $j=1, 2, ..., n,$ define $y^{j}:B_{\delta}(\bar{p})\to\mathbb{R}$ by $y^{j}=\pi^{j}\circ\varphi$, where $\pi^{j}:T_{\bar{p}}\mathbb{M} \to\mathbb{R}$ is the projection defined by $\pi^{j}\left(a_{1}, ..., a_{j}, ..., a_{n}\right)=a_{j}, $
for all $ \left(a_{1}, ..., a_{j}, ..., a_{n}\right)\in T_{\bar{p}}\mathbb{M}$. 
Then,   $\left( B_{\delta}(\bar{p}),\varphi, y^{j} \right)$  is a   local coordinate system at $\bar{p}$. Let  $\{ \partial/ \partial y^{j} \}$   the associated  corespondent  natural  basis to $\left( B_{\delta}(\bar{p}),\varphi, y^{j} \right)$.  Since  $\gamma_{p}(t) \in  B_{\delta}(\bar{p})$  for all $t\in [0, 1]$ and  $Y_{p}\left(\gamma_{p}(t)\right)\in T_{\gamma_{p}(t)}\mathbb{M}$ we can write
$$
Y_{p}\left(\gamma_{p}(t)\right)=\underset{j}\sum\, Y^{j}_{p}(t)\dfrac{\partial}{\partial y^{j}}\vert_{\gamma_{p}(t)}\quad \forall ~ t\in[0,1],
$$
where each coordinate function  $Y^{j}_{p}:[0,1]\to\mathbb{R}$ is  continuously differentiable, for all  $j=1, 2, ..., n$. For simplicity, we set $y^{j}_{p}:=y^{j}\circ\gamma_{p}(.)$ for each $j=1, 2, ..., n$.  Thus \eqref{derivadanula} is equivalent to  the ordinary  differential equation
$$
\dfrac{dY^{k}_{p}}{dt}+\underset{i,j}{\sum}\,\Gamma^{k}_{i,j}(\gamma_{p})\dfrac{dy^{i}_{p}}{dt}Y^{j}_{p}=0,\quad k=1, 2, ..., n,
$$
where $\Gamma^{k}_{i,j}$ is the Christoffel symbols of the connection $\nabla$, see \cite[pag. 29]{Sakai1996}.  Hence,  last equality implies that  $\{Y^{k}_{p}~: k=1, 2, ..., n\}$ is the unique solution of the following the system of  $p-$parameter linear differential equations
\[ 
 \left\{ \begin{array}{ll}
            \dfrac{dY^{k}_{p}}{dt}=-\displaystyle{\sum^{n}_{j}\,a_{k,j}Y^{j}_{p}},\quad k=1, 2, ..., n,\\
            \displaystyle{\sum_{j}\, Y^{j}_{p}(0)\dfrac{\partial}{\partial y^{j}}\vert_{\gamma_{p}(0)}=u},                                                                    
\end{array} \right. 
\]
 where, for  $(k,j)$ with $k,j=1, ..., n,$  the continuous  function  $a_{k,j}:[0,1]\times B_{\delta}(\bar p)\to\mathbb{R}$ is given by  
$$
a_{k,j}(t,p)=\sum^{n}_{i=1}\,\Gamma^{k}_{i,j}\left(\gamma_{p}(t)\right)\dfrac{dy^{i}_{p}(t)}{dt}.
$$
Thus, from result about continuity on  parameters for  differential equations  (see,  for example, \cite[Theorem~{10.7.1}, pag. 353]{DieudonnE1969}),  the solution $\{Y^{k}_{.}(\cdot )\}$ is continuous on $[0,1]\times B_{\delta}(\bar{p})$ and equivalently, $Y_{.}(\gamma_{.}(\cdot))$ is continuous on $[0,1]\times B_{\delta}(\bar{p})$. Furthermore, we have $ F(p)=P_{\bar p p}u=Y_{p}(\gamma_{p}(1))$, for any $p\in B_{\delta}(\bar{p}).$ Therefore,  $F$ is continuous on $B_{\delta}(\bar{p})$ and  the  is completed.
\qed

Next,  we present an immediate consequence of  Lemma~\ref{le:ContTranpPar}. 
\begin{corollary}\label{cor:clongLiinvertido} 
Let    $\bar{p} \in  \mathbb{M}$, $0<\delta<i_{\bar{p}}$ and  $u\in T_{\bar{p}}\mathbb{M}$.  Then, if the vector field   $Z: B_{\delta}(\bar{p})  \to T\mathbb{M}$ is  continuous  at $\bar{p}$, then  the mapping $G: B_{\delta}(\bar{p}) \to  T_{\bar{p}}\mathbb{M} $ defined by $G(p):=P_{p\bar{p} }Z(p)$ is also  continuous at $\bar{p}$.
\end{corollary}
{\it Proof} Since the parallel transport is a isometry, it follows from  the definition of  vector field $G$ that 
$$
\| G(p)-G(\bar p)\|=\| Z(p)-P_{\bar p p}Z(\bar p)\|.
$$
Taking into account that $Z$ is continuous at $\bar p$ and  $P_{\bar p \bar p}=I_{\bar p}$, we conclude from Lemma \ref{le:ContTranpPar} that 
$$
\underset{p\to\bar p}{\lim}\| Z(p)-P_{\bar p p}Z(\bar p)\|=0.
$$
Therefore, the  desired result  follows by simple combination of the two last equalities.
\qed

The next result ensure us that, if $\nabla X(p_{*})$ is nonsingular then there exist a neighborhood of $p_*$ which $\nabla X$ also is nonsingular.  Besides,  in this neighborhood $\nabla X^{-1}$ is bounded.
\begin{lemma}\label{le:NonSing}
 Assume that $\nabla X$ is  continuous at $p_{*}$. Then, there holds
 \begin{equation}\label{Hdepvaiazero}
\underset{p\to p_{*}}{\lim}\left \| P_{pp_{*}}\nabla X(p)P_{p_{*}p}-\nabla X(p_{*}) \right \|=0.
\end{equation}
Moreover, if  $\nabla X(p_{*})$ is  nonsingular,  then there exists $0<\bar{\delta}<i_{p_{*}}$ such that $B_{\bar{\delta}}(p_{*})\subset \Omega$  and  for each $p\in B_{\bar{\delta}}(p_{*})$ there hold:
\begin{itemize}
\item[i)] $\nabla X(p)$ is nonsingular;
\item[ii)] $\left \| \nabla X(p)^{-1}\right \|\leq 2\left \|\nabla X(p_{*})^{-1}\right \|$.
\end{itemize}
\end{lemma}
{\it Proof} Let $0<\delta <  i_{p_{*}}$ such that $B_{\delta}(p_{*})\subset\Omega$. For each $u\in T_{p_{*}}\mathbb{M}$,  define  $Z: B_{\delta}(p_{*}) \to T\mathbb{M}$  by
$$
Z(p)=\nabla X(p)P_{p_{*}p}u.
$$
Applying   Lemma~\ref{le:ContTranpPar}  we  conclude   that  $P_{p_{*}p}u$ is continuous on $B_{\delta}(p_{*})$. Thus,  considering that
$\nabla X$ is continuous, we obtain that  $Z$ is also continuous on $B_{\delta}(p_{*})$. Hence, using  Corollary~\ref{cor:clongLiinvertido},  we conclude  that the mapping $F:B_{\delta}(p_{*})\to T_{p_{*}}\mathbb{M}$ defined by
$$
F(p)=P_{pp_{*}}Z(p), 
$$ 
is also continuous at $p_{*}$. Taking  into account  that  $P_{p_{*}p_{*}}=I_{p_{*}}$ and the  definitions of  the mappings $F$ and $Z$, we conclude  $\underset{p\to p_{*}}{\lim}F(p)=\nabla X(p_{*})u$. Now, define the mapping
$$
B_{\delta}(p_{*})\owns p\mapsto \left[P_{pp_{*}}\nabla X(p)P_{p_{*}p}-\nabla X(p_{*})\right]\in  {\cal L}(T_{p_{*}}\mathbb{M},T_{p_{*}}\mathbb{M}), 
$$
where   ${\cal L}(T_{p_{*}}\mathbb{M},T_{p_{*}}\mathbb{M}) $ denotes  the space consisting of all  linear operator  from $T_{p_{*}}\mathbb{M}$ to  $T_{p_{*}}\mathbb{M}$. Since $\underset{p\to p_{*}}{\lim}F(p)=\nabla X(p_{*})u$,  for each $u\in T_{p_{*}}\mathbb{M}$, thus the  definitions of  $F$ implies 
$$
\underset{p\to p_{*}}{\lim}\left[ P_{pp_{*}}\nabla X(p)P_{p_{*}p}-\nabla X(p_{*}) \right]u=0,\quad u\in T_{p_{*}}\mathbb{M}.
$$
Owing to the fact that $T_{p_{*}}\mathbb{M}$  is finite dimensional and $[P_{pp_{*}}\nabla X(p)P_{p_{*}p}-\nabla X(p_{*})] \in {\cal L}(T_{p_{*}}\mathbb{M},T_{p_{*}}\mathbb{M})$,  for each $p\in B_{\delta}(p_{*})$,  the latter equality   implies  that the  equality   \eqref{Hdepvaiazero} holds. Now,  we proceed with the proof of the first item. The equality \eqref{Hdepvaiazero} implies that there exists $0<\bar\delta<\delta$ such that
$$
\left\| P_{pp_{*}}\nabla X(p)P_{p_{*}p}-\nabla X(p_{*}) \right\|\leq\dfrac{1}{2\left\|\nabla X(p_{*})^{-1}\right\|},\qquad  \forall ~p\in B_{\bar\delta}(p_{*}).
$$
Thus,  from  the last inequality  and the property  of norm of operator defined in ${\cal L}(T_{p_{*}}\mathbb{M},T_{p_{*}}\mathbb{M})$, for all $p\in B_{\bar\delta}(p_{*})$ we obtain
\begin{equation}\label{Bpnonsingular}
\left\|\nabla X(p_{*})^{-1}P_{pp_{*}}\nabla X(p)P_{p_{*}p}-I_{p_{*}}\right\| \leq
\| \nabla X(p_{*})^{-1}\|\left\|P_{pp_{*}}\nabla X(p)P_{p_{*}p}-\nabla X(p_{*}) \right\|
\leq \dfrac{1}{2}.
\end{equation}
Hence from Lemma \ref{le:Banach} we conclude that $\nabla X(p_{*})^{-1}P_{pp_{*}}\nabla X(p)P_{p_{*}p}$ is nonsingular operator for each $p\in B_{\bar{\delta}}(p_{*})$. Due to   $\nabla X(p_{*})$  and the parallel transport are nonsingular  we obtain $\nabla X(p)$ is also nonsingular for each $p\in B_{\bar{\delta}}(p_{*})$ and the proof  of the first item is concluded. To prove item $ii$ we first note that, from  \eqref{Bpnonsingular} and  Lemma \ref{le:Banach}, for all $p\in B_{\bar{\delta}}(p_{*})$ follows that
$$
\left\| \left[\nabla X(p_{*})^{-1}P_{pp_{*}}\nabla X(p)P_{p_{*}p}\right]^{-1} \right\|\leq\dfrac{1}{1-\left\|\nabla X(p_{*})^{-1}P_{pp_{*}}\nabla X(p)P_{p_{*}p}-I_{p_{*}} \right\|}. 
$$
Since the parallel transport is an isometry, combining \eqref{Bpnonsingular} with     the latter inequality we obtain that 
$$
\left\|  \nabla X(p)^{-1} P_{p_{*}p} \nabla X(p_{*})  \right\|\leq 2.   
$$
Thus, using  the  properties of the norm and again that the parallel transport is an isometry, the  last  inequality  implies that, for all $p\in B_{\bar{\delta}}(p_{*})$ we have 
$$
\left\|  \nabla X(p)^{-1} \right\|\leq \left\|   \nabla X(p)^{-1} P_{p_{*}p} \nabla X(p_{*})\right\| \left\|\nabla X(p_{*})^{-1}\right\|  \leq 2\left\|\nabla X(p_{*})^{-1}\right\|,   
$$ 
which is the desired  inequality in the second item.   Thus, the proof of the lemma  is  concluded.
\qed

Lemma~\ref{le:NonSing} establishes  non-singularity of $\nabla X$ in a neighborhood of $p_{*}$.  It ensures us that there exists a neighborhood of $p_{*}$ which Newton's iterate \eqref{eq:NM} is well defined. But it not guarantee  that  it belongs to this neighborhood. In the next lemma, we will establish it. For stating  the next result, we first  define the {\it Newton's iterate mapping}   $N_{X}:B_{\bar{\delta}}(p_{*}) \to  \mathbb{M} \nonumber $  by
\begin{equation} \label{Newton's iterate mapping}
N_{X}(p):=\exp_{p}(-\nabla X(p)^{-1}X(p)), 
\end{equation}
where $\bar{\delta}$  is  given by Lemma~\ref{le:NonSing}.
\begin{lemma}\label{L5}
Assume  that $\nabla X$ is continuous at $p_{*}$ and $\nabla X(p_{*})$ is  nonsingular.  Then, there holds
$$
\lim_{p\to p_{*}}\frac{d(N_{X}(p),p_{*})}{d(p,p_{*})}=0.
$$
\end{lemma}
{\it Proof} Let   $\bar{\delta} $  be given by Lemma \ref{le:NonSing} and   $p\in B_{\bar{\delta}}(p_{*})$. Some algebraic manipulations  show that
\begin{multline*}
 \nabla X(p)^{-1}X(p)+\exp^{-1}_pp^*= 
\nabla X(p)^{-1} \big{[ }X(p)- P_{p_{*}p}X(p_{*})-\\
P_{p_{*}p}\nabla X(p_{*})\exp^{-1}_{p^*}p+ \left[  P_{p_{*}p}\nabla X(p_{*})-\nabla X(p)P_{p_{*}p}\right]\exp^{-1}_{p^*}p\big{]}.
\end{multline*}
Define $r(p):= [X(p)- P_{p_{*}p}X(p_{*})-P_{p_{*}p}\nabla X(p_{*})\exp^{-1}_{p^*}p]/d(p, p^*)$, for $p\in B_{\bar{\delta}}(p_{*})$.  From   \eqref{diferenciabilidade} we have  ${\lim_{p\to p_{*}}}r(p)=0$. Thus, using the above equality, the definition of $r$,  $d(p, p^*)=\| \exp^{-1}_{p^*}p\|$ and some  properties of the norm,  we  conclude  that
$$
\left\|\nabla X(p)^{-1}X(p)+\exp^{-1}_pp^* \right\|\leq\left\|\nabla X(p)^{-1}\right \| \left[ \left\| r(p) \right\|+
\left\| P_{p_{*}p}\nabla X(p_{*})-\nabla X(p)P_{p_{*}p}\right\| \right]d(p, p^*).
$$
Since $p\in   B_{\bar{\delta}}({p_{*}})$ and the parallel transport is an isometry,  thus  item ii of Lemma \ref{le:NonSing}  implies  that
\begin{equation} \label{eq:hel}
\left\|\nabla X(p)^{-1}X(p)+\exp^{-1}_pp^* \right\| \leq  2\left\|\nabla X(p_{*})^{-1}\right\|  \left[ \left\| r(p) \right\|+ \left\| P_{pp_{*}}\nabla X(p)P_{p_{*}p}-\nabla X(p_{*})\right\| \right]d(p, p^*).
\end{equation}
Due to \eqref{Hdepvaiazero} and  $ \underset{p\to p_{*}}{\lim}r(p)=0$,  the right hand side of last inequality goes to zero, as $p$ goes to $p^*$. Thus we  can shrink $ \bar\delta$,  if necessary,  in order to obtain  that
$$
\left\|  \nabla X(p)^{-1}X(p)+\exp^{-1}_pp^* \right\|\leq i_{p_{*}},\qquad  \forall ~p\in   B_{ \bar\delta}(p_{*}).
$$
Hence, from definition of Newton's iterate mapping $N_{X}$ in \eqref{Newton's iterate mapping} and definition  of $K_{p_{*}} $ we have
$$
d(N_{X}(p),p_{*})\leq  K_{p_{*}}\left\| - \nabla X(p)^{-1}X(p)-\exp^{-1}_pp^* \right\|,\qquad  \forall ~p\in   B_{ \bar\delta}(p_{*}).
$$
Therefore, combining \eqref{eq:hel} with  the last  inequality  we conclude, for all $p\in B_{ \bar\delta}(p_{*})$ the following
$$
\dfrac{d(N_{X}(p),p_{*})}{d(p, p_{*})}\leq 2 K_{p_{*}} \left\|\nabla X(p_{*})^{-1}\right\|  [\parallel r(p)\parallel+\parallel P_{pp_{*}}\nabla X(p)P_{p_{*}p}-\nabla X(p_{*})\parallel].
$$
Letting $p$ goes to $p{*}$ in the last inequality, taking into account  \eqref{Hdepvaiazero} and that  ${\lim_{p\to p_{*}}}r(p)=0$,  the desired result follows.
\qed

Now,  we are ready to establishes  our  main result, its proof will be a combination of the two previous lemmas.
\begin{theorem}
Let $\Omega\subset\mathbb{M}$ be a open set, $X:\Omega\to T\mathbb{M}$ a differentiable vector field and $p_{*}\in\Omega$. Suppose that $p_{*}$ is a singularity of $X$, $\nabla X$ continuous at $p_{*}$ and $\nabla X(p_{*})$ is  nonsingular. Then,  there exist $\bar{\delta}>0$ such that,  for  all  $p_{0}\in B_{\bar{\delta}}(p_{*})$,  the Newton's sequence 
\begin{equation} \label{eq:NM2}
p_{k+1}=\exp_{p_{k}}(-\nabla X(p_{k})^{-1}X(p_{k})),  \qquad k=0, 1, \ldots, 
\end{equation}
is well defined, contained in $B_{\bar{\delta}}(p_{*})$ and converges superlinearly to $p_{*}$.
\end{theorem}
{\it Proof}
Let  $ \bar{\delta} $ be  given by  Lemma~\ref{le:NonSing}.  From Lemma~\ref{L5}   we  can shrink $ \bar\delta$,  if necessary,
to conclude that  
\begin{equation}\label{convergencia} 
d(N_{X}(p), p_{*})<\dfrac{1}{2}d(p, p_{*}),\qquad \forall ~p\in B_{\bar{\delta}}(p_{*}).
\end{equation}
Thus, $N_{X}(p)\in B_{\bar{\delta}}(p_{*})$, for all $p\in B_{\bar{\delta}}(p_{*})$.  Note that \eqref{Newton's iterate mapping}  together   \eqref{eq:NM2}  implies that  $\{p_k\}$  satisfies 
\begin{equation} \label{NFS}
p_{k+1}=N_X(p_k),\qquad k=0,1,\ldots \,,
\end{equation}
which is indeed an equivalent definition of this sequence. Since $N_{X}(p)\in B_{\bar{\delta}}(p_{*})$, for all $p\in B_{\bar{\delta}}(p_{*})$,   it follows from \eqref{NFS} and  Lemma \ref{le:NonSing} item $i)$  that, for all $p_{0}\in B_{\bar{\delta}}(p_{*})$   the Newton sequence $\{p_{k}\}$ is well defined and  contained in $B_{\bar{\delta}}(p_{*})$. Moreover, using   (\ref{convergencia}) and \eqref{NFS} we obtain  that 
$$
d(p_{k+1}, p_{*})<\dfrac{1}{2}d(p_{k}, p_{*}),\qquad  k=0,\,1, \ldots.
$$
The latter inequality implies that  $\{p_{k}\}$ converges to $p_{*}$. Thus,   combining   \eqref{NFS}  with     Lemma~\ref{L5}, we conclude that 
$$
\lim_{k\to +\infty}\frac{d(p_{k+1},p_{*})}{d(p_{k},p_{*})}=0. 
$$
Therefore, $\{p_{k}\}$  convergence superlinearly to $p_{*}$ and the proof is concluded.
\qed

Before concluding  this section, we present  an  important property of the parallel transport.  When $\mathbb{M}$ is a complete and finite dimensional Riemannian manifold, the Lemma~\ref{le:ContTranpPar} allow us  obtain the continuity of the vector field $B_{\delta}(\bar p)\ni p \mapsto P_{p{\bar p}}Z(p)$, where $\delta\leq i_{\bar p}$;  see Corollary~\ref{cor:clongLiinvertido}. On the other hand, when $\mathbb{M}$ is the Hadamard manifold, i.e., a complete simply connected Riemannian manifold of non-positive sectional curvature, we know that $i_{\bar p}=+\infty$,   $\phi:=exp_{\bar p}^{-1}:\mathbb{M}\to T_{\bar p}\mathbb{M}$ is a diffeomorfism and  $(\mathbb{M},\phi)$ is a global chart for $\mathbb{M}$. Therefore,  following the same idea of the proof  Corollary~\ref{cor:clongLiinvertido}  we can prove the following  generalization of \cite [Lemma 2.4, item (iv)]{LiLopezMartin-Marquez2009}:
 \begin{corollary}
Let $\mathbb{M}$ be Hadamard manifold.  If  $Z:\mathbb{M}\to T\mathbb{M}$ is a continuous vector field, then the mapping $ \mathbb{M}\ni(p,q) \mapsto P_{p,q}Z(p)$ is continuous.
\end{corollary}
{\it Proof} First, from \cite [Lemma 2.4, item (iv)]{LiLopezMartin-Marquez2009} we have ${\lim_{q\to\bar q}}P_{\bar p q}Z(\bar p)=P_{\bar p\bar q}Z(\bar p)$. Then,  to prove the result  is enough show that the following equality holds
\begin{equation}\label{limiteequivalente}
\underset{(p,q)\to(\bar p,\bar q)}{\lim}\|P_{pq}Z(p)-P_{\bar p q}Z(\bar p)\|=0.
\end{equation}
 Indeed, after some algebraic manipulations  and  considering that  $P_{\bar p q}=P_{\bar q q}P_{\bar p \bar q}$, $P_{q\bar q}P_{\bar q q}=I_{\bar q}$, $P_{p\bar q}P_{\bar q p}=I_{\bar q}$ and the parallel transport is a isometry, we have
$$
\|P_{pq}Z(p)-P_{\bar p q}Z(\bar p)\|=\|Z(p)-P_{\bar p p}Z(\bar p)\|.
$$
Since,  \cite [Lemma 2.4, item (iv)]{LiLopezMartin-Marquez2009} implies that $\underset{p\to\bar p}{\lim}\|Z(p)-P_{\bar p p}Z(\bar p)\|=0$, we have from the last equality that  \eqref{limiteequivalente} holds and the proof is concluded. 
\qed

\section{Final Remarks}

In this paper,   under  non-singularity of the covariant derivative of the vector field at its  zero  and  without any additional conditions on this derivative, we establish the  superlinear local convergence of Newton's method  on a finite dimensional  Riemannian manifold.  It is worth to pointed out that  we have  assumed that the Riemannian manifold is finite dimensional to establishes  this result, at least two times, namely,  in Lemma~\ref{le:ContTranpPar}  and Lemma~\ref{le:NonSing}. On the other hand, we know that   Newton's method  converges  superlinearly  on infinite dimensional Banach space under non-singularity of the derivative of the operator  at its  zero.  It would be interesting to extend this results to  infinite dimensional Riemannian manifolds.




\begin{thebibliography}{10}
\providecommand{\url}[1]{{#1}}
\providecommand{\urlprefix}{URL }
\expandafter\ifx\csname urlstyle\endcsname\relax
  \providecommand{\doi}[1]{DOI~\discretionary{}{}{}#1}\else
  \providecommand{\doi}{DOI~\discretionary{}{}{}\begingroup
  \urlstyle{rm}\Url}\fi

\bibitem{Smith1994}
Smith, S.T.: Optimization techniques on {R}iemannian manifolds.
\newblock In: Hamiltonian and gradient flows, algorithms and control,
  \emph{Fields Inst. Commun.}, vol.~3, pp. 113--136. Amer. Math. Soc.,
  Providence, RI (1994)

\bibitem{Absil2009}
Absil, P.A., Mahony, R., Sepulchre, R.: Optimization algorithms on matrix
  manifolds.
\newblock Princeton University Press, Princeton, NJ (2008)

\bibitem{Luenberger1972}
Luenberger, D.G.: The gradient projection method along geodesics.
\newblock Management Sci. \textbf{18}, 620--631 (1972)

\bibitem{Gabay1982}
Gabay, D.: Minimizing a differentiable function over a differential manifold.
\newblock J. Optim. Theory Appl. \textbf{37}(2), 177--219 (1982)


\bibitem{Rapcsak1989}
Rapcs{\'a}k, T.: Minimum problems on differentiable manifolds.
\newblock Optimization \textbf{20}(1), 3--13 (1989)


\bibitem{Udriste1994}
Udri{\c{s}}te, C.: Convex functions and optimization methods on {R}iemannian
  manifolds, \emph{Mathematics and its Applications}, vol. 297.
\newblock Kluwer Academic Publishers Group, Dordrecht (1994)


\bibitem{Absil2014}
Absil, P.A., Amodei, L., Meyer, G.: Two {N}ewton methods on the manifold of
  fixed-rank matrices endowed with {R}iemannian quotient geometries.
\newblock Comput. Statist. \textbf{29}(3-4), 569--590 (2014)


\bibitem{Huang2015}
Huang, W., Gallivan, K.A., Absil, P.A.: A {B}royden class of quasi-{N}ewton
  methods for {R}iemannian optimization.
\newblock SIAM J. Optim. \textbf{25}(3), 1660--1685 (2015)


\bibitem{LiLopezMartin-Marquez2009}
Li, C., L{\'o}pez, G., Mart{\'{\i}}n-M{\'a}rquez, V.: Monotone vector fields
  and the proximal point algorithm on {H}adamard manifolds.
\newblock J. Lond. Math. Soc. (2) \textbf{79}(3), 663--683 (2009)


\bibitem{LiWang2008}
Li, C., Wang, J.: Newton's method for sections on {R}iemannian manifolds:
  generalized covariant {$\alpha$}-theory.
\newblock J. Complexity \textbf{24}(3), 423--451 (2008)


\bibitem{WangLiLopezYao2015}
Wang, J., Li, C., Lopez, G., Yao, J.C.: Convergence analysis of inexact
  proximal point algorithms on {H}adamard manifolds.
\newblock J. Global Optim. \textbf{61}(3), 553--573 (2015)


\bibitem{WangLiWangYao2015}
Wang, X., Li, C., Wang, J., Yao, J.C.: Linear convergence of subgradient
  algorithm for convex feasibility on {R}iemannian manifolds.
\newblock SIAM J. Optim. \textbf{25}(4), 2334--2358 (2015)


\bibitem{HuLiYang2016}
Hu, Y., Li, C., Yang, X.: On convergence rates of linearized proximal
  algorithms for convex composite optimization with applications.
\newblock SIAM J. Optim. \textbf{26}(2), 1207--1235 (2016)


\bibitem{Ring2012}
Ring, W., Wirth, B.: Optimization methods on {R}iemannian manifolds and their
  application to shape space.
\newblock SIAM J. Optim. \textbf{22}(2), 596--627 (2012)


\bibitem{Argyros2014}
Argyros, I.K., Magre{\~n}{\'a}n, {\'A}.A.: Extending the applicability of
  {G}auss-{N}ewton method for convex composite optimization on {R}iemannian
  manifolds.
\newblock Appl. Math. Comput. \textbf{249}, 453--467 (2014)


\bibitem{Manton2015}
Manton, J.H.: A framework for generalising the {N}ewton method and other
  iterative methods from {E}uclidean space to manifolds.
\newblock Numer. Math. \textbf{129}(1), 91--125 (2015)


\bibitem{Bittencourt2016}
Bittencourt, T., Ferreira, O.: Kantorovich’s theorem on newton’s method
  under majorant condition in riemannian manifolds.
\newblock Journal of Global Optimization pp. 1--25 (2016)

\bibitem{De2016}
de~Carvalho~Bento, G., da~Cruz~Neto, J.X., Oliveira, P.R.: A new approach to
  the proximal point method: convergence on general {R}iemannian manifolds.
\newblock J. Optim. Theory Appl. \textbf{168}(3), 743--755 (2016)


\bibitem{Adler2002}
Adler, R.L., Dedieu, J.P., Margulies, J.Y., Martens, M., Shub, M.: Newton's
  method on {R}iemannian manifolds and a geometric model for the human spine.
\newblock IMA J. Numer. Anal. \textbf{22}(3), 359--390 (2002)

\bibitem{Nash1956}
Nash, J.: The imbedding problem for {R}iemannian manifolds.
\newblock Ann. of Math. (2) \textbf{63}, 20--63 (1956)

\bibitem{Moser1961}
Moser, J.: A new technique for the construction of solutions of nonlinear
  differential equations.
\newblock Proc. Nat. Acad. Sci. U.S.A. \textbf{47}, 1824--1831 (1961)

\bibitem{Owren2000}
Owren, B., Welfert, B.: The {N}ewton iteration on {L}ie groups.
\newblock BIT \textbf{40}(1), 121--145 (2000)


\bibitem{Li2009}
Li, C., Wang, J.H., Dedieu, J.P.: Smale's point estimate theory for {N}ewton's
  method on {L}ie groups.
\newblock J. Complexity \textbf{25}(2), 128--151 (2009)


\bibitem{Argyros2009}
Argyros, I.K., Hilout, S.: Newton's method for approximating zeros of vector
  fields on {R}iemannian manifolds.
\newblock J. Appl. Math. Comput. \textbf{29}(1-2), 417--427 (2009)


\bibitem{Dedieu2000}
Dedieu, J.P., Shub, M.: Multihomogeneous {N}ewton methods.
\newblock Math. Comp. \textbf{69}(231), 1071--1098 (electronic) (2000)


\bibitem{Miller2005}
Miller, S.A., Malick, J.: Newton methods for nonsmooth convex minimization:
  connections among $\mathcal{U}$-agrangian, {R}iemannian {N}ewton and {SQP}
  methods.
\newblock Math. Program. \textbf{104}(2-3, Ser. B), 609--633 (2005)


\bibitem{Schulz2014}
Schulz, V.H.: A {R}iemannian view on shape optimization.
\newblock Found. Comput. Math. \textbf{14}(3), 483--501 (2014)


\bibitem{Wang2011}
Wang, J.H.: Convergence of {N}ewton's method for sections on {R}iemannian
  manifolds.
\newblock J. Optim. Theory Appl. \textbf{148}(1), 125--145 (2011)

\bibitem{Wang2009}
Wang, J.H., Li, C.: Kantorovich's theorems for {N}ewton's method for mappings
  and optimization problems on {L}ie groups.
\newblock IMA J. Numer. Anal. \textbf{31}(1), 322--347 (2011)


\bibitem{Wang2012}
Wang, J.H., Yao, J.C., Li, C.: Gauss-{N}ewton method for convex composite
  optimizations on {R}iemannian manifolds.
\newblock J. Global Optim. \textbf{53}(1), 5--28 (2012)


\bibitem{Zhang2010}
Zhang, L.H.: Riemannian {N}ewton method for the multivariate eigenvalue
  problem.
\newblock SIAM J. Matrix Anal. Appl. \textbf{31}(5), 2972--2996 (2010)


\bibitem{Ferreira2012}
Ferreira, O.P., Silva, R.C.M.: Local convergence of {N}ewton's method under a
  majorant condition in {R}iemannian manifolds.
\newblock IMA J. Numer. Anal. \textbf{32}(4), 1696--1713 (2012)


\bibitem{Ferreira2002}
Ferreira, O.P., Svaiter, B.F.: Kantorovich's theorem on {N}ewton's method in
  {R}iemannian manifolds.
\newblock J. Complexity \textbf{18}(1), 304--329 (2002)


\bibitem{Dedieu2003}
Dedieu, J.P., Priouret, P., Malajovich, G.: Newton's method on {R}iemannian
  manifolds: convariant alpha theory.
\newblock IMA J. Numer. Anal. \textbf{23}(3), 395--419 (2003)


\bibitem{Jiang2007}
Jiang, D., Moore, J.B., Ji, H.: Self-concordant functions for optimization on
  smooth manifolds.
\newblock J. Global Optim. \textbf{38}(3), 437--457 (2007)


\bibitem{Li2006}
Li, C., Wang, J.: Newton's method on {R}iemannian manifolds: {S}male's point
  estimate theory under the {$\gamma$}-condition.
\newblock IMA J. Numer. Anal. \textbf{26}(2), 228--251 (2006)


\bibitem{Ortega1990}
Ortega, J.M.: Numerical analysis. {A} second course.
\newblock Academic Press, New York-London (1972)


\bibitem{doCarmo1992}
do~Carmo, M.P.: Riemannian geometry.
\newblock Mathematics: Theory \& Applications. Birkh\"auser Boston, Inc.,
  Boston, MA (1992)


\bibitem{Sakai1996}
Sakai, T.: Riemannian geometry, \emph{Translations of Mathematical Monographs},
  vol. 149.
\newblock American Mathematical Society, Providence, RI (1996)


\bibitem{DieudonnE1969}
Dieudonn{\'e}, J.: Foundations of modern analysis.
\newblock Academic Press, New York-London (1969)


\end{thebibliography}

\end{document}